\documentclass[aps,prl,preprint,groupedaddress]{revtex4-1}
\usepackage{graphicx}
\usepackage{dcolumn}
\usepackage{bm}
\newcommand{\rmi}{\textrm{i}}
\newcommand{\rmd}{\textrm{d}}
\newcommand{\rme}{\textrm{e}}
\newtheorem{corollary}{\hskip\parindent\bf Corollary}
\newtheorem{theorem}{\hskip\parindent\bf Theorem}
\begin{document}

\preprint{}

\title[Codimension-two Bifurcations Induce Hysteresis Behavior and Multistabilities]{Codimension-two Bifurcations Induce Hysteresis Behavior and Multistabilities in  Delay-coupled Kuramoto Oscillators}

\author{Ben Niu}
\email{niubenhit@163.com}
\affiliation{Department of  Mathematics, Harbin Institute of Technology (Weihai),\\  Weihai 264209, P.R. China}
\date{\today}

\begin{abstract}
Hysteresis phenomena and multistability play crucial roles in the dynamics of  coupled oscillators, which are now interpreted from the point of view of codimension-two bifurcations. On the Ott-Antonsen's manifold, complete bifurcation sets of delay-coupled Kuramoto model are derived regarding coupling strength and delay as  bifurcation parameters.  It is  rigorously proved that the system  must undergo Bautin bifurcations for some critical values, thus there always exists saddle-node bifurcation of periodic solutions inducing hysteresis loop. With the aid of center manifold reduction method and the Matlab Package DDE-Biftool, the location of Bautin and double Hopf points and detailed dynamics are theoretically determined.  We find that, near these critical points, at most four coherent states (two of which are stable) and a stable incoherent state may coexist, and that the system  undergoes Neimark-Sacker  bifurcation of periodic solutions. Finally, the clear scenarios about the synchronous transition in delayed Kuramoto model are depicted.
\end{abstract}

\maketitle

\section {Introduction}

The Kuramoto model was  established to investigate the phenomenon of
collective synchronization of coupled oscillators with slightly different natural frequencies
\cite{ref1,ref2,ref3,ref4}, which has been widely observed in physics, chemistry and biology
\cite{article10,article10-1,article10-2,article10-3,article10-4,article10-5}. Because of the time lag that signal transmits or that the receiver   processes the signal, introducing time delay is natural
and necessary in many situations \cite{ar9,ar91,ar901,ar912,ott1,ott2,wai,arkd,niunon}.
The Kuramoto  model with time delay is of the form
\begin{equation}\label{model}
\dot
\theta_i=\omega_i+\frac{k}{N}\sum_{j=1}^{N}\sin[\theta_j(t-\tau)-\theta_i(t)],i=1,2,\ldots,N\end{equation}
where $\theta_i(t)\in[0,2\pi)$ represents the phase of the $i$th oscillator at time $t$. $\omega_i,i=1,2,\ldots,N$ are  natural frequencies  drawn from
  density $g(\omega), \omega\in(-\infty,+\infty)$,  and  positive $k$ is
the coupling strength.  In the thermodynamic limit $N\rightarrow\infty$,   define a distribution density $f(\omega,\theta,t)$ characterizing the state of the oscillators' system at time $t$ in frequency $\omega$ and phase $\theta$. Then the complex-valued mean field $r(t)$ is define by
\begin{equation}\label{rrrr}
r(t)=\int_{-\infty}^{+\infty}\int_0^{2\pi}f(\omega,\theta,t)\rme^{\rmi\theta}\rmd\theta \rmd\omega
\end{equation}
describing the degree to which the oscillators are bunched in phase, where $\rmi=\sqrt{-1}$. Write the continuity equation
\begin{equation}\label{fp}
\frac\partial{\partial t}f+\frac\partial{\partial \theta}\{[\omega+\frac k{2\rmi}(\rme^{-\rmi\theta}r(t-\tau)-\rme^{\rmi\theta}r^\ast(t-\tau))]
f\}=0\end{equation}
with $r^\ast$ standing for the complex conjugate. Usually we call $r=0$ ($f$ is uniform distribution) the incoherent state, $|r|=1$ ($f$ is Dirac distribution) the completely synchronized state, and $0<|r|<1$ the partially synchronized state (coherence, for short).

Hopf bifurcation research in Kuramoto model is an efficient way to obtain the transition between the incoherence and coherence \cite{ott1,niuphyD}. Coherent states bifurcating from the incoherent state can be modeled by Hopf bifurcations on some low dimensional manifold, such as the widely used Ott-Antonsen's manifold \cite{ott1}. The direction of Hopf bifurcation, subcritical case or supercritical case,  then determines different situations of the synchronization transition.   

In \cite{ar9} the hysteresis loop and subcritical bifurcations are observed in the delay coupled Kuramoto oscillators. Here, when  a hysteresis loop is mentioned, we mean that coherent states and incoherent states coexist in the Kuramoto model when the parameter $k$ is less than the Hopf bifurcation value. In \cite{niuphyD}, the authors have interpreted the appearance of subcritical Hopf bifurcations in the way of normal form analysis. However, a clear boundary between the supercritical and subcritical bifurcations (a degenerated case) has not been theoretically studied yet. In the viewpoint of bifurcation analysis, this may be involved with the Bautin bifurcation with codimension-two (i.e., generalized Hopf bifurcation)\cite{guo,bautin,bautin1,bautin2}. This is a degenerated case we mainly considered in the current paper. Another degenerated case occurs when two Hopf bifurcation coexist, i.e., the double Hopf bifurcation, which is also codimention two and rarely investigated in the Kuramoto model before. It is well-known that double Hopf bifurcation usually provides a system with oscillations on a 2-torus or 3-torus \cite{guo,homs,kuz} through the Neimark-Sacker bifurcation of periodic solutions. To our best knowledge, codimension-two bifurcation (including Bautin and double Hopf bifurcations) approach to dynamical analysis  is brand new to investigate delayed Kuramoto model.

Motivated by such considerations, in this paper, we study the Bautin bifurcation  and double Hopf bifurcation on the Ott-Antonsen's manifold  \cite{ott1} to reveal some delicate dynamics for
delay coupled system (\ref{model}).
The rest part of this paper is organized as follows:  we first restate the OA manifold reduction method with respect to system (\ref{fp}), and analyze the characteristic equation of the incoherence. Then Bautin bifurcation and double Hopf bifurcation are analyzed with the aid of center manifold reduction method, respectively.  Some illustrations are given with the help of the Matlab Package DDE-Biftool, hence a clear bifurcation set is given in the $\tau-k$ plane. The results are also applied to a system of delay-coupled Hindmarsh-Rose neurons.
Finally a conclusion part completes this paper.

\section{OA Manifold Reduction and Stability of the Incoherence}
 For the readers' convenience, we first restate the main results about OA manifold reduction of (\ref{model}) by \cite{ott1}.
Restrict (\ref{fp}) on the OA manifold
\[
\mathcal{M}_{OA}=\left\{f(\theta,\omega, t):f=\frac{g(\omega)}{2\pi}\left\{1+\left[\sum_{n=1}^\infty\alpha^n(\omega,t)\rme^{\rmi n\theta}+c.c.\right]\right\}\right\}\]
with c.c. the complex conjugate of the formal terms and $\alpha^n(\omega,t)$ the Fourier coefficients.
Substituting the Fourier series  of $f$ into  (\ref{fp}) and after comparing the coefficients of the  same harmonic terms, a reduced equation is obtained
\begin{equation}\label{alpha}
\dot \alpha(\omega,t)=-\rmi \omega\alpha(\omega,t)+\frac k 2 r^\ast(t-\tau)-\frac k 2 r(t-\tau) \alpha^2(\omega,t)
\end{equation}
Following the OA ansatz \cite{ott1} and applying Cauchy's residue theorem to Eq.(\ref{rrrr}),   we have \begin{equation}\label{rt}r(t)=\int_{-\infty}^{+\infty}g(\omega)\alpha^\ast(\omega,t)\rmd\omega = \alpha^\ast(\omega_0-\rmi\Delta,t)\end{equation}
provided that
$g(\omega)$ is  chosen to be the Lorentzian distribution
$
g(\omega)=\frac{\Delta}{\pi[(\omega-\omega_0)^2+\Delta^2]},~-\infty<\omega<+\infty$
with the spreading  width $\Delta>0$ and the median value $\omega_0>0$.

Putting $\omega=\omega_0-\rmi\Delta$ in Eq.(\ref{alpha}) yields a reduced functional differential equation \cite{Halefde}
\begin{equation}\label{model1}
\dot r(t)=-(\rmi\omega_0+\Delta)r(t)+\frac k 2 r(t-\tau)-\frac k 2 r^\ast(t-\tau)r^2(t)
\end{equation}

To give the transition from incoherence to coherence, we need to investigate
the
characteristic equation around the incoherence $r=0$
\begin{equation}\label{CElaplacedelta}
\lambda=-(\rmi\omega_0+\Delta)+\frac k 2 \rme^{-\lambda
\tau}\end{equation}

Since we are about to study codimension-two bifurcations, regarding $k$ and $\tau$ as bifurcation parameters, we know Eq.(\ref{model1}) undergoes local bifurcation at $r=0$ if (\ref{CElaplacedelta}) has roots with zero real part for some $(k,\tau)$. Noticing the assumption $\omega_0\neq0$, we let
$\lambda=\rmi\beta$ with $\beta\neq0$ be a root of (\ref{CElaplacedelta}), then
\begin{equation}\label{cossin}\begin{array}{rcl}
\Delta &=&\frac{1}{2} k \cos \tau \beta \\
-\beta -\omega_0 &=&\frac{1}{2} k \sin \tau \beta \end{array}
\end{equation}
which yields \[\beta ^2+2   \omega_0
\beta+\Delta ^2+\omega_0^2 -\frac{k^2}{4}=0\]
Obviously,  $\beta$ is solved by
\[
\beta=\beta_\pm=-\omega_0\pm\sqrt{\frac {k^2}{4}-\Delta^2}
\]
if $k>2\Delta$ holds.
Furthermore, from (\ref{cossin}), two sequences of critical values $\tau_j^\pm(k)$ are defined by
\begin{equation}\label{tau}
\tau_j^+(k)=\left\{\begin{array}{ll}\frac{\arcsin \frac {-2\beta_+-2\omega_0}k-2j\pi}{\beta_+},&k^2\leq4(\Delta^2+\omega_0^2)\\
\frac{\arcsin \frac {-2\beta_+-2\omega_0}k+2(j+1)\pi}{\beta_+},&k^2>4(\Delta^2+\omega_0^2)\end{array}j=0,1,2,\cdots\right.
\end{equation}
and
\begin{equation}\label{tau2}
\tau_j^-(k)=\frac{\arcsin \frac {-2\beta_--2\omega_0}k-2(j+1)\pi}{\beta_-},j=0,1,2,\cdots
\end{equation}
This means that (\ref{CElaplacedelta}) has a root $\rmi\beta_+$ (or $\rmi\beta_-$), if $(\tau,k)=(\bar\tau,\bar k)\in\{(\tau,k)|k>2\Delta,\tau=\tau_j^+{(k)}~(\textrm{or}~\tau_j^-(k))\}$. Usually, purely imaginary roots of characteristic equation mean Hopf bifurcation or Bautin bifurcation. The bifurcating periodic solution $r(t)\approx r_0\rme^{\rmi\beta t}$ with small amplitude $r_0$ of (\ref{model1}) corresponds to a coherent state of (\ref{model}). Obtaining precise results requires the normal forms near the critical points.

\section{Bautin bifurcation}
 When $k=\bar k>2\Delta$, the characteristic equation (\ref{CElaplacedelta}) has a purely imaginary root if $\tau=\bar\tau\in\{\tau_j^\pm(\bar k),j=0,1,2,\cdots\}$. In order to obtain the bifurcation results, we need calculate the normal form by using the center manifold reduction method \cite{Hassard,faria,guo} basing on the formal adjoint theory \cite{Halefde}. It is worthy mentioning that the method of multiple time scales can be also used to obtained normal forms in delay equations \cite{naybook,Yums,dfmms,Nay}.  The two approaches lead to the same normal forms \cite{ding}, thus we use the center manifold reduction approach here.
Normalizing time by   $\tau$, we rewrite Eq.(\ref{model1}) as
\begin{equation}\label{model2}\begin{array}{rl}
\dot r(t)=-\tau(\rmi\omega_0+\Delta)r(t)+\frac k 2 \tau r(t-1)-\frac k 2 \tau r^\ast(t-1)r^2(t)
\end{array}
\end{equation}
with a characteristic equation at $r
=0$
\begin{equation}\label{CE}
\lambda=-\tau(\rmi\omega_0+\Delta)+\tau\frac  k 2 \rme^{-\lambda
}\end{equation}
Clearly, $\rmi\beta$ is a root of (\ref{CElaplacedelta}) if and only if $\rmi\tau\beta$ is a   root of (\ref{CE}).
Slightly perturbing $k=\bar k+\epsilon$ and $\tau=\bar\tau+\delta$, we have the equivalent form of (\ref{model2})
\begin{equation}\label{model3}\begin{array}{rl}
\dot r=-(\bar\tau+\delta)(\rmi\omega_0+\Delta)r+\frac {(\bar k+\epsilon)(\bar\tau+\delta)} 2  r(t-1)-\frac {(\bar k+\epsilon)(\bar\tau+\delta)} 2  r^\ast(t-1)r^2\end{array}
\end{equation}
If $(\epsilon,\delta)=(0,0)$, i.e., $(\tau,k)=(\bar\tau,\bar k)$, $\rmi\bar\tau\beta$ is a  root of (\ref{CE}).

For any $\varphi
\in C([-1,0],\mathbf{C})$ with $\mathbf{C}$ the  space of complex numbers, we define a linear operator
\[
L_{\epsilon,\delta}(\varphi)=
                                 \int_{-1}^0\rmd \eta(\theta,\epsilon,\delta)\varphi(\theta) \\
                            \]
with
\[
 \eta(\theta,\epsilon,\delta)=\left\{\begin{array}{ll}-(\bar\tau+\delta)(\rmi\omega_0+\Delta)+\frac {(\bar k+\epsilon)} 2 (\bar\tau+\delta),&\theta=0\\
 \frac {(\bar k+\epsilon)} 2 (\bar\tau+\delta),&\theta\in(-1,0)\\
0,&\theta=-1
 \end{array}\right.
\]
Meanwhile, for $\varphi\in C^1([-1,0],\mathbf{C})$, we define
\[
A(\epsilon,\delta)\varphi=\left\{\begin{array}{ll}\frac{\rmd\varphi(\theta)}{\rmd\theta},&\theta\in[-1,0)\\
L_{\epsilon,\delta}(\varphi),&\theta=0\end{array}\right.
\]
\[
F(\epsilon,\delta,\varphi)=-\frac {(\bar k+\epsilon)(\bar\tau+\delta) } 2                                                                                \varphi^\ast(-1)\varphi^2(0)
\]
and
\[
R(\epsilon,\delta)\varphi=\left\{\begin{array}{ll}0,&\theta\in[-1,0)\\
F(\epsilon,\delta,\varphi),&\theta=0\end{array}\right.
\]
Then system (\ref{model3}) can be transformed into an abstract ordinary differential equation
\begin{equation}\label{reducedeq}
\dot r_t=A(\epsilon,\delta)r_t+R(\epsilon,\delta)r_t
\end{equation}
with $r_t=                          r(t+\theta),\theta\in[-1,0]$.

For $\psi\in C^1([-1,0],{\mathbf{C}})$, define the adjoint operator of $A(0,0)$ by
\[
\hat A\psi =\left\{\begin{array}{ll}-\frac{\rmd\psi(s)}{\rmd s},&s\in(0,1]\\
\int_{-1}^0\psi(-s)\rmd \eta^\ast(\theta,0,0),& s=0\end{array}\right.
\]
For $\varphi\in C([-1,0],\mathbf{C})$ and $\psi\in C([0,1],{\mathbf{C}})$, define the bilinear form
\[
<\psi,\phi>=\psi^\ast(0)\varphi(0)-\int_{\theta=-1}^0\int_{\zeta=0}^\theta\psi^\ast(\zeta-\theta)\rmd\eta(\theta,0,0)\varphi(\zeta)\rmd\zeta
\]
We know that $\rmi\bar\tau\beta$ is an eigenvalue of $A(0,0)$ and $-\rmi\bar\tau\beta$ is an eigenvalue of $\hat A$. Suppose $q(\theta)$  and $\hat q(s)$ are the corresponding eigenvectors, i.e.,
\[\begin{array}{l}
A(0,0)q=\rmi\bar\tau\beta q,\hat A \hat q=-\rmi\bar\tau\beta \hat q
\end{array}\]
Letting
$q(\theta)=
\rme^{\rmi\beta\bar\tau\theta}$, $\hat q(s)= \frac{1}{1+\frac{1}{2} \bar k \bar \tau  \rme^{\rmi \beta  \bar\tau }} \rme^{\rmi\beta\bar\tau s}$, we have $<\hat q,q>=1$.

Due to the classical results in \cite{Halefde}, for $|(\epsilon,\delta)|$ sufficiently small, we use $z(t)$ as  complex coordinate on the center manifold in direction $q$,   thus  $z(t)=<\hat q,r_t>$.
Decomposing \begin{equation}\label{decomrt}r_t(\theta)=z(t)q(\theta) +W^{(\epsilon,\delta)}(z,z^\ast,\theta)\end{equation}
with
$W^{(\epsilon,\delta)}(z,z^\ast,\theta)=\sum_{i+j\geq2}W_{ij}^{(\epsilon,\delta)}(\theta)\frac{1}{i!j!}z^i{z^\ast}^j$.
Denote $W^{(0,0)}$ by $W$ and $W_{ij}^{(0,0)}$ by $W_{ij}$ for simplicity. We can calculate
\begin{equation}\label{modelcenter}\begin{array}{rl}\dot z(t)&=<\hat q,\dot r_t>\\&=<\hat q,A(0,0)r_t+[A(\epsilon,\delta)-A(0,0)]r_t+R(\epsilon,\delta)r_t>\\&=<\hat A\hat q,r_t>+<\hat q,[A(\epsilon,\delta)-A(0,0)]r_t>+<\hat q,R(\epsilon,\delta)r_t>\\&=
\rmi\bar\tau\beta z(t)+<\hat q,[A(\epsilon,\delta)-A(0,0)]r_t>+\hat q^\ast(0)F(\epsilon,\delta,r_t)\end{array}\end{equation}
Rewrite (\ref{modelcenter}) as
\begin{equation}\label{zeq}
\dot z(t)=\rmi\bar\tau\beta z(t)+<\hat q,[A(\epsilon,\delta)-A(0,0)]r_t>+g(z,z^\ast)
\end{equation}
then we Taylor expand $g(z,z^\ast)=\hat q^\ast(0)F(\epsilon,\delta,r_t):=\sum_{i+j\geq2}g_{ij} \frac{1}{i!j!}z^i{z^\ast}^j$.
According to the results about normal form with imaginary roots \cite{guo,Wiggins2}, we have the
 normal form is
\[
 \dot z=\rmi\bar\tau\beta z+l_0(\epsilon,\delta)z +\frac12l_1(\epsilon,\delta)z^2z^\ast+\frac1{12}l_2(\epsilon,\delta)z^3{z^{\ast}}^{2}+h.o.t.
\]
 with
\[l_0(\epsilon,\delta)= \frac{1}{1+\frac{1}{2}\bar k \bar\tau  \rme^{-\rmi \beta \bar \tau }}\left[\frac{1}{2}   \bar \tau  \rme^{-\rmi \beta  \bar \tau }\epsilon+ \left(-\Delta +\frac{1}{2}\bar  k \rme^{-\rmi \beta \bar  \tau }-\rmi \omega _0\right)\delta\right]\]
Obviously, $(Re\frac {\partial {l_0}}{\partial \epsilon},Re\frac {\partial {l_0}}{\partial \delta})\neq0$ yields the satisfaction of transversality condition. The first and second Lyapunov coefficients at $\epsilon=\delta=0$ are
$
Rel_1(0,0)=\frac {Img_{20} g_{11} } {\bar\tau|\beta|} + {Reg_{21}}
$ and
\[\begin{array}{rl}
&12Rel_2(0,0)=Re{g_{32}}\\&~~+\frac1{\bar\tau|\beta|}Im\left\{g_{20}\bar g_{31}-g_{11}(4g_{31}+3\bar g_{22})-\frac13g_{02}(g_{40}+\bar g_{13})-g_{30}g_{12}\right\}\\
&~~+\frac{1}{\bar\tau^2\beta^2}Re\left\{g_{20}\left[\bar g_{11}(3g_{12}-\bar g_{30})+g_{02}\left(\bar g_{12}-\frac13g_{30}\right)+\bar g_{02}g_{03}\right]\right\}\\
&~~+\frac{1}{\bar\tau^2\beta^2}Re\left\{g_{11}\left[\bar g_{02}\left(\frac53\bar g_{30}+3g_{12}\right)+\frac13g_{02}\bar g_{03}-4g_{11}g_{30}\right]\right\}\\
&~~+\frac{3}{\bar\tau^2\beta^2}Im\{g_{20}g_{11}\}Im\{g_{21}\}+\frac1{\bar\tau^3|\beta|^3}Im\{
g_{11}\bar g_{02}[\bar g_{20}^2-3\bar g_{20}g_{11}-4g_{11}^2\}\\
&~~+\frac1{\bar\tau^3|\beta|^3}Im\left\{
g_{11} g_{20}\right\}[3Re\{g_{11}g_{20}\}-2|g_{02}|^2]
 \end{array}\]

 Due to the lack of second order terms in $F(\epsilon,\delta,\varphi)$, simply we have \[g_{20}=g_{11}=g_{02}=0,~g_{21}=-\frac{\bar k\bar\tau}{2} \frac{\rme^{\rmi \beta  \bar\tau }}{1+\frac{1}{2} \bar k\bar \tau  \rme^{-\rmi \beta  \bar\tau }}\]
\[g_{30}=g_{12}=g_{03}=g_{04}=0,~g_{40}=-\frac{\bar k\bar\tau}{2}\frac1{{1+\frac{1}{2} \bar k\bar \tau  \rme^{-\rmi \beta  \bar\tau }}}(w_{02}^\ast(-1))\]
\[g_{31}=-\frac{\bar k\bar\tau}{2}\frac1{{1+\frac{1}{2} \bar k\bar \tau  \rme^{-\rmi \beta  \bar\tau }}}(2\rme^{-\rmi\beta\bar\tau}w_{20}(0)+w_{11}^\ast(-1))\]
\[g_{22}=-\frac{\bar k\bar\tau}{2}\frac1{{1+\frac{1}{2} \bar k\bar \tau  \rme^{-\rmi \beta  \bar\tau }}}(2\rme^{-\rmi\beta\bar\tau}w_{11}(0)+w_{20}^\ast(-1))\]
\[g_{13}=-\frac{\bar k\bar\tau}{2}\frac1{{1+\frac{1}{2} \bar k\bar \tau  \rme^{-\rmi \beta  \bar\tau }}}(2\rme^{-\rmi\beta\bar\tau}w_{02}(0))\]and
\[g_{32}=-\frac{\bar k\bar\tau}{2}\frac1{{1+\frac{1}{2} \bar k\bar \tau  \rme^{-\rmi \beta  \bar\tau }}}(w_{21}^\ast(-1)+2w_{20}^\ast(-1)w_{20}(0)+2w_{02}^\ast(-1)w_{02}(0)+2\rme^{-\rmi\beta\bar\tau}w_{21}(0))\]
Plugging  (\ref{decomrt}) and (\ref{modelcenter}) into (\ref{reducedeq}), we have
\begin{equation}\label{WWW}
\dot W=\left\{\begin{array}{l}AW-\hat q^\ast(0)F(0,0,r_t)q(\theta),-1\leq\theta<0\\AW-\hat q^\ast(0)F(0,0,r_t) q(\theta)+F(0,0,r_t),\theta=0\end{array}\right.
\end{equation}
In fact $\dot W=W_{z}\dot z+W_{z^\ast}\dot z^\ast$, then balancing the coefficients  $z$ and $z^\ast$ in (\ref{WWW}), we have some differential equations and initial conditions, by solving which we obtain the unknown $w_{20}, w_{02}, w_{11}$ and $w_{21}$ above. Here we omit the tedious expressions (these calculations can be found in \cite{guo,bautin1,bautin2}) and only give the final expressions of $l_1(0,0)$ and $l_2(0,0)$ by running a computer program,
\[l_1(0,0) =-\frac{\bar k\bar\tau}{2} Re\frac{\rme^{\rmi \beta  \bar\tau }}{1+\frac{1}{2} \bar k\bar \tau  \rme^{-\rmi \beta  \bar\tau }},l_2(0,0)= - 3\bar k\bar\tau  Re \frac{1}{1+\frac{1}{2}\bar k \bar\tau  \rme^{-\rmi \beta \bar \tau }}\]

  Using the fundamental results in \cite{homs,kuz,Wiggins2}, we have
 \begin{theorem}
At the critical point $(\bar k,\bar\tau)$, if Re$l_1(0,0)\neq0$, then system (\ref{model1}) undergoes a Hopf bifurcation, which is supercritical if Re$l_1(0,0)<0$ and subcritical if Re$l_1(0,0)>0$. If Re$l_1(0,0)=0$, and Re$l_2(0,0)\neq0$, then system (\ref{model1}) undergoes a  Bautin bifurcation.
\end{theorem}

\begin{figure}[htbp]
  \centering \includegraphics[width=0.6\textwidth]{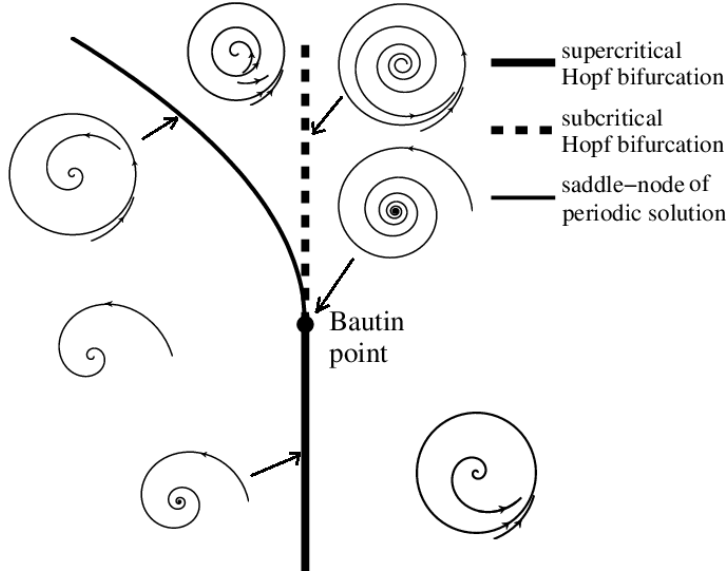}
  \caption{Bifurcation sets around a generic Bautin point with Re$l_2(0,0)<0$.}\label{bautinfig}
\end{figure}

In fact Sgn Re $l_1(0,0)=$ Sgn $-\Delta-2\bar\tau\Delta^2+\frac {\bar\tau} 4{\bar k}^2$, thus we have
 \begin{corollary}\label{corol}For any $\Delta>0,\omega_0>0$, system (\ref{model1}) always undergoes  Bautin bifurcations at some $(\bar\tau,\bar k)$. There always exists a curve of saddle-node bifurcations of periodic orbits originating from each Bautin point, which gives birth to the hysteresis loop in the delay-coupled Kuramoto model (\ref{model}). A generic result is shown in Figure \ref{bautinfig}.\end{corollary}
{\bf Proof:} In the $\tau-k$ plane, we know Hopf bifurcation appears for all $ k>2\Delta$.
 The curve $-\Delta-2\tau\Delta^2+\frac {\tau} 4{ k}^2=0$ admits $k=\sqrt{{4\Delta}/{\tau}+8\Delta^2}$, which is a monotone decreasing function of $\tau$, and approximates to $ k=2\sqrt{2}\Delta>2\Delta$ as $\tau\rightarrow\infty$, and to $+\infty$ as $\tau\rightarrow0$. Thus we know there must exist some intersecting points between the Hopf curves and the critical curve  $-\Delta-2\tau\Delta^2+\frac {\tau} 4{ k}^2=0$, i.e., Bautin points.
Moreover, we have Re$l_2(0,0)=$Re$[1+\bar \tau(\rmi\beta+\rmi\omega_0+\Delta)]^{-1}>0$ from (\ref{CE}).  The rest results are direct applications of  the classical results in \cite{homs}.

\section{Double Hopf Bifurcation}
Usually, two Hopf bifurcation curves on the $\tau-k$ plane intersect, leaving the system with some double Hopf points \cite{guo}.

 Still denoting these double Hopf points by $(\bar\tau,\bar k)$, we set $(\tau,k)=(\bar\tau,\bar k)+(\epsilon,\delta)$ to proceed bifurcation analysis. We assume in this case that $\rmi\bar\tau\alpha$ and $\rmi\bar\tau\beta$ are eigenvalues of $A(0,0)$. Hence $-\rmi\bar\tau\alpha$ and  $-\rmi\bar\tau\beta$ are eigenvalues of $\hat A$. Suppose $q_1(\theta),q_2(\theta)$ and $\hat q_1(s),\hat q_2(s)$ are the corresponding eigenvectors, i.e.,
\[
A(0,0)q_1=\rmi\bar\tau\alpha q,\hat A \hat q_1=-\rmi\bar\tau\alpha \hat q,
A(0,0)q_2=\rmi\bar\tau\beta q,\hat A \hat q_2=-\rmi\bar\tau\beta \hat q
\]
Let
$q_1(\theta)=
\rme^{\rmi\alpha\bar\tau\theta}$, $q_2(\theta)= \rme^{\rmi\beta\bar\tau\theta}$ and $\hat q_1(s)= D_1^\ast\rme^{\rmi\alpha\bar\tau s}$, $\hat q_2(s)= D_2^\ast\rme^{\rmi\beta\bar\tau s}$, $
    D_1=\frac{1}{1+\frac{1}{2} \bar k \bar\tau  \rme^{-\rmi \alpha  \bar\tau }},
    D_2=\frac{1}{1+\frac{1}{2} \bar k \bar\tau  \rme^{-\rmi \beta  \bar\tau }}$, then we have $<\hat q_1,q_1>=<\hat q_2,q_2>=1$ and $<\hat q_1,q_2>=<\hat q_2,q_1>=0$.

Using $z_1(t)$ and $z_2(t)$ as  complex coordinates on the center manifold for small $|(\epsilon,\delta)|$, we have $z_1(t)=<\hat q_1,r_t>$, $z_2(t)=<\hat q_2,r_t>$.
Then
\[\begin{array}{rl}\dot z_1(t)&=<\hat q_1,\dot r_t>\\&=<\hat q_1, A(\epsilon,\delta)r_t+R(\epsilon,\delta)r_t>\\&=<\hat A\hat q_1,r_t>+<\hat q_1, [A(\epsilon,\delta)-A(0,0)])r_t>+<\hat q_1,R(\epsilon,\delta)r_t>\\&=
\rmi\bar\tau\alpha z_1(t)+<\hat q_1, [A(\epsilon,\delta)-A(0,0)])r_t>+\hat q_1^\ast(0)F(\epsilon,\delta,r_t)\end{array}\]
\[\begin{array}{rl}\dot z_2(t)&=<\hat q_2,\dot r_t>\\&=<\hat q_2, A(\epsilon,\delta)r_t+R(\epsilon,\delta)r_t>\\&=<\hat A\hat q_2,r_t>+<\hat q_2, [A(\epsilon,\delta)-A(0,0)])r_t>+<\hat q_2,R(\epsilon,\delta)r_t>\\&=
\rmi\bar\tau\beta z_2(t)+<\hat q_2, [A(\epsilon,\delta)-A(0,0)])r_t>+\hat q_2^\ast(0)F(\epsilon,\delta,r_t)\end{array}\]
Letting $r_t(\theta)=z_1(t)q_1(\theta)+z_2(t)q_2(\theta)+W(z,z^\ast,\theta)$, and following the same method given in \cite{guo},
 if the nonresonant condition
 \begin{equation}\label{nonre}|\alpha|\neq 3|\beta|, 3|\alpha|\neq |\beta|\end{equation}is satisfied, the third order normal form near a double Hopf point is derived
\[\begin{array}{lll}
      \dot z_1 &=&a_{11}\epsilon z_1+a_{12}\delta z_1+c_{11}z_1^2 z_1^\ast+c_{12}z_1z_2  z_2^\ast\\
    \dot z_2 &=&a_{21}\epsilon z_2+a_{22}\delta z_2+c_{21}z_1  z_1^\ast z_2+c_{22}z_2^2  z_2^\ast\end{array}
\]
where
\[\begin{array}{l}
   a_{11} = \frac{1}{2} D_1 \bar\tau  \rme^{-\rmi \alpha  \bar\tau },
   a_{12} = D_1 \left(-\Delta +\frac{1}{2}\bar k \rme^{-\rmi \alpha  \bar\tau }-\rmi \omega _0\right)\\
   c_{11} = -\frac{1}{2} D_1 \bar k \bar\tau  \rme^{\rmi \alpha  \bar\tau },
  c_{12} =-\frac{1}{2}D_1\bar k \bar\tau  \rme^{\rmi \beta  \bar\tau } \\
    a_{21} =\frac{1}{2} D_2 \bar\tau  \rme^{-\rmi \beta \bar \tau },
   a_{22} = D_2 \left(-\Delta +\frac{1}{2}\bar k \rme^{-\rmi \beta  \bar\tau }-\rmi \omega _0\right)\\
   c_{21} = -\frac{1}{2} D_2\bar k \bar\tau  \rme^{\rmi \alpha \bar \tau },
  c_{22} =-\frac{1}{2} D_2\bar k \bar\tau  \rme^{\rmi \beta \bar \tau }\end{array}
\]

After   rescaling    $\epsilon_1=\textrm{Sign}(\textrm{Re}c_{11})$,
$\epsilon_2=\textrm{Sign}(\textrm{Re}c_{22})$,  $r_1=|z_1|,$ $r_2=|z_2|$,  $t \rightarrow t \epsilon_1$ we have the amplitude equation
\begin{equation}\label{NF23polarreducehh}\begin{array}{lll}
       \dot r_1&=&r_1(c_1+r_1^2+b_0r_2^2)\\ \dot r_2&=& r_2(c_2+c_0r_1^2+d_0r_2^2)
\end{array}
\end{equation}
with
\[
\begin{array}{l}
 c_1= \epsilon_1\textrm{Re}a_{11}\epsilon+\epsilon_1\textrm{Re}a_{12}\delta,
 c_2 = \epsilon_1\textrm{Re}a_{21}\epsilon
+\epsilon_1\textrm{Re}a_{22}\delta \\
  b_0 = \frac{\epsilon_1\epsilon_2\textrm{Re}c_{12}}{\textrm{Re}c_{22}},
 c_0 = \frac{\textrm{Re}c_{21}}{\textrm{Re}c_{11}},
 d_0 = \epsilon_1\epsilon_2
 \end{array}
\]

Applying the results in \cite{homs}, Eq.(\ref{NF23polarreducehh})   has twelve distinct types of unfoldings,
distinguished by the signs of $b_0,~c_0,~d_0$ and $d_0-b_0c_0$. In the coming section, we will give a numerical example to show the detailed bifurcation sets near a double Hopf point.

\section{Numerical Experiments}

In this section, some illustrations are given to support the theoretical results obtained about Bautin and double Hopf bifurcations in (\ref{model2}). Meanwhile the Kuramoto model (\ref{model}) is also simulated. The existence of multistabilities is observed in a delay-coupled system of Hindmarsh-Rose neurons.

\subsection{Simulations near the Bautin bifurcation}
\begin{figure}[htbp]
  \centering a)\includegraphics[width=0.47\textwidth,height=0.31\textwidth]{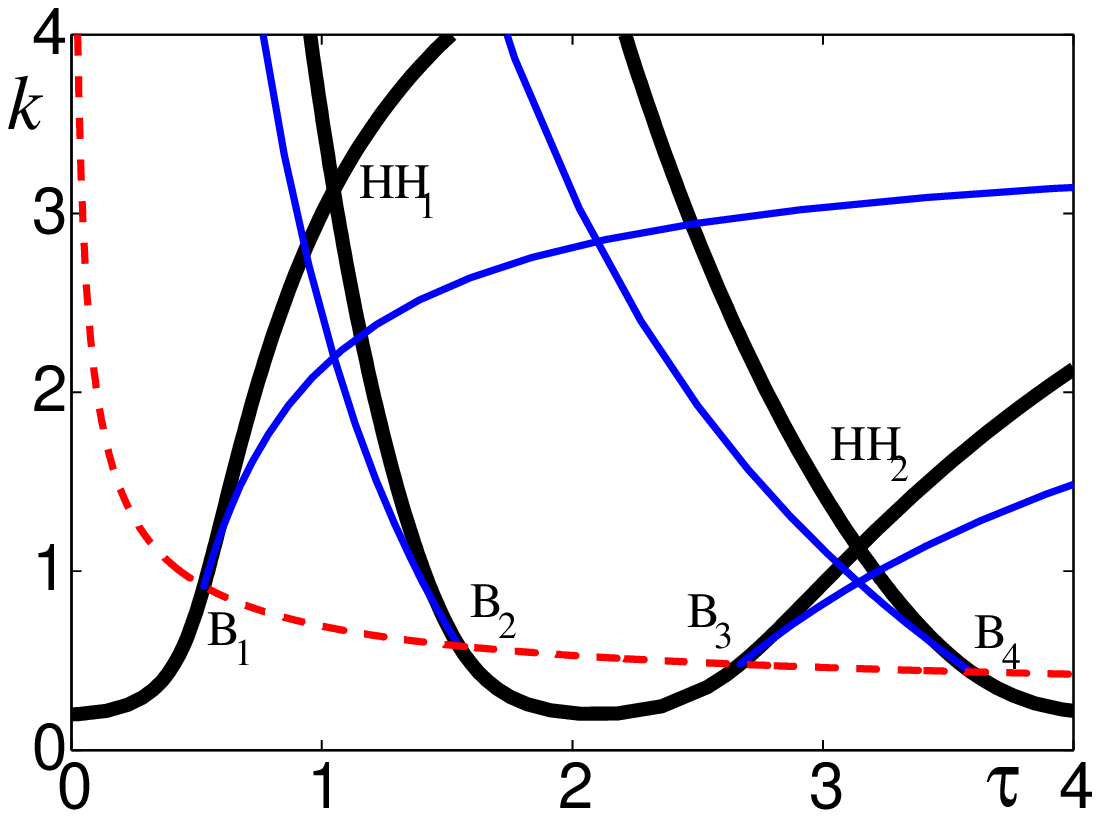}b)\includegraphics[width=0.47\textwidth,height=0.31\textwidth]{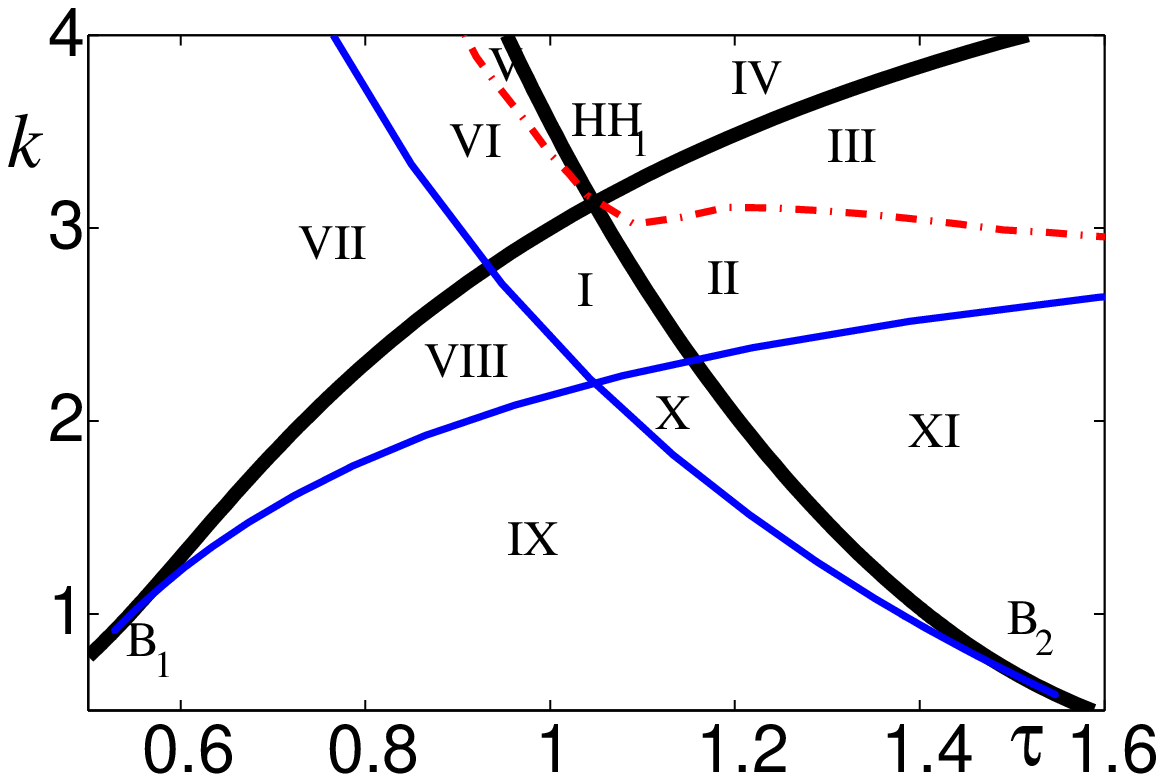}
  \\c)\includegraphics[width=0.47\textwidth]{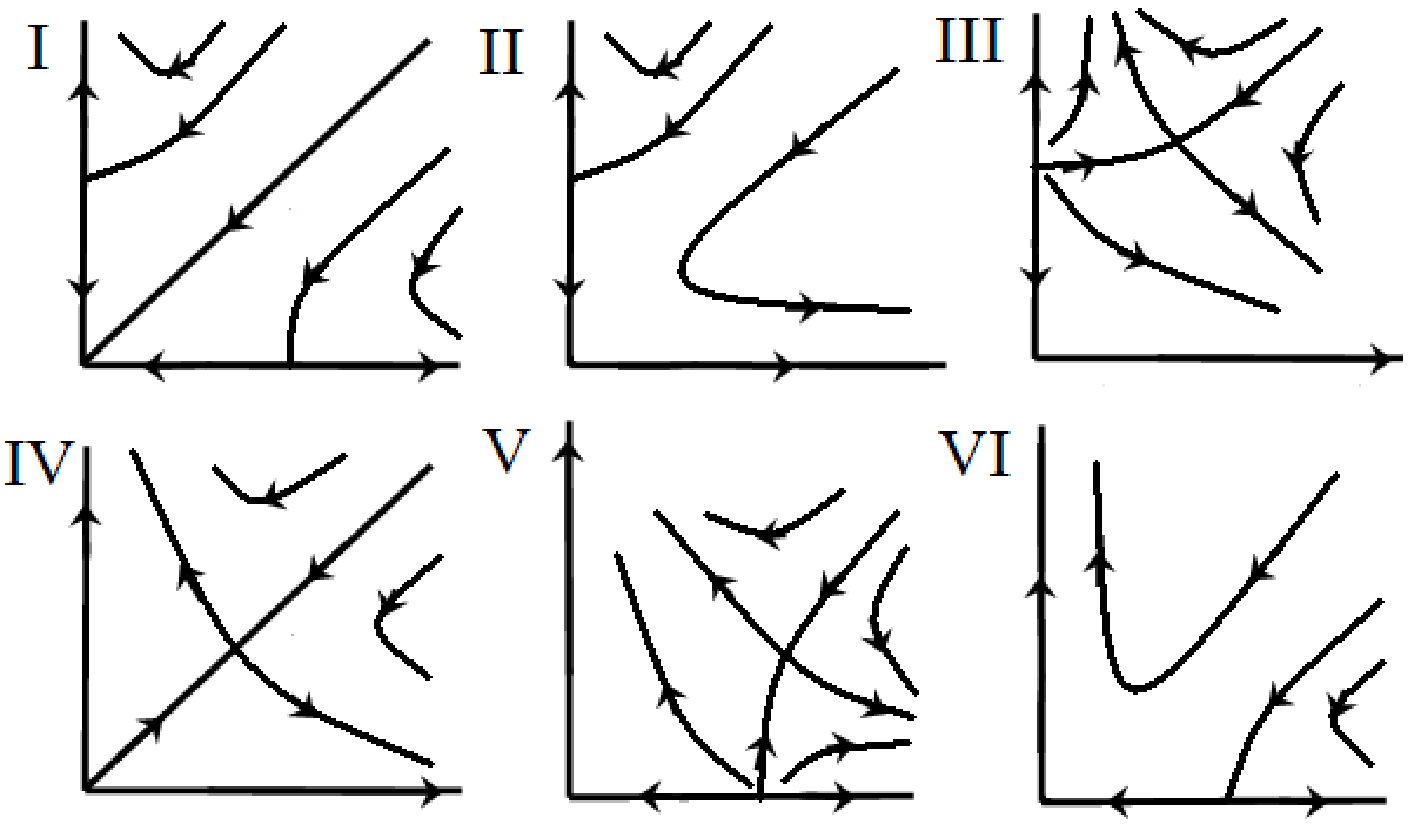}d)\includegraphics[width=0.47\textwidth]{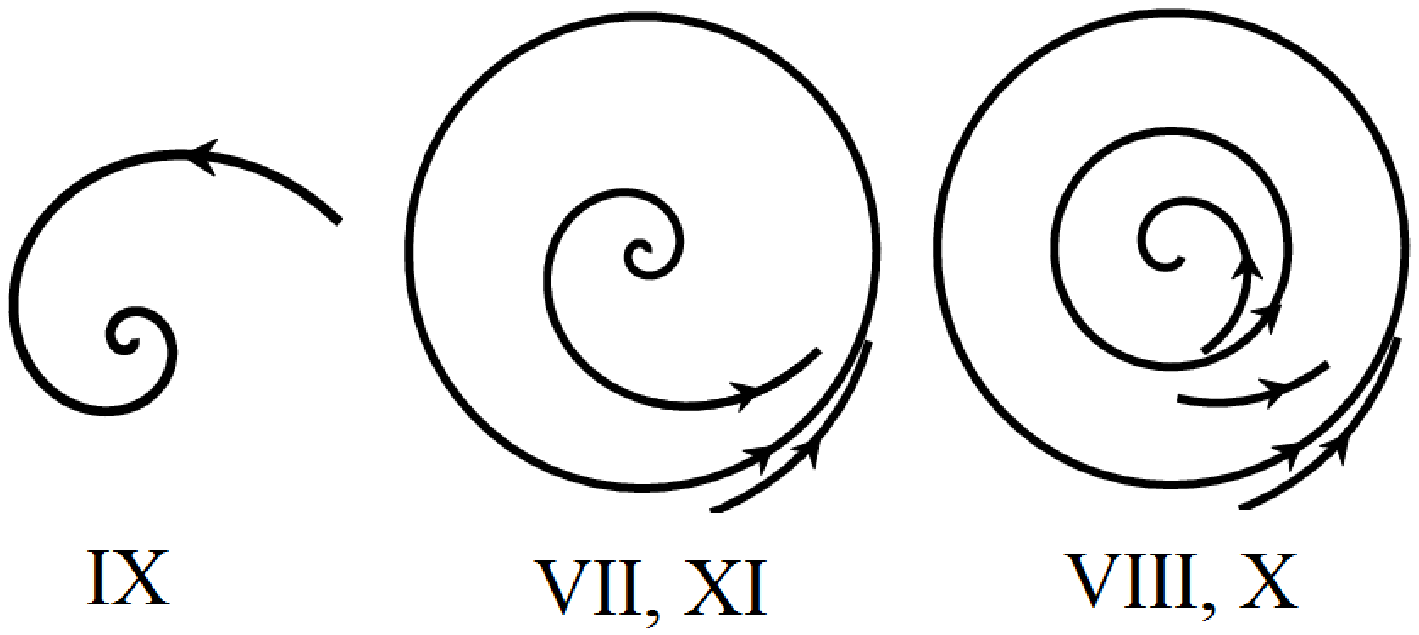}
  \caption{(color online) a) Hopf bifurcation curves (black, thick) of Eq.(\ref{model1}) are shown in the $k-\tau$ plane for $\omega_0=3$ and $\Delta=0.1$. HH$_1-$HH$_2$ are double Hopf points. B$_1-$B$_4$ are Bautin points. Blue (thin) curves  stand for saddle-node bifurcation of periodic orbits. Red (dashed) curve stands for the line $-\Delta-2\tau\Delta^2+\frac \tau4k^2=0$, above which Hopf bifurcations are subcritical and below which Hopf bifurcations are supercritical.  b) figure a) is zoomed in the region $[0.5,1.6]\times [0.5,4]$ and   the red dash-dotted curves stand  for the Neimark-Sacker  bifurcations. c) local phase portraits for parameters in regions I-VI of b) in polar coordinates $r_1$ and $r_2$ of Eq.(\ref{NF23polarreducehh}). d) phase portraits for parameters in regions VII-XI.}\label{Tkkkk}
\end{figure}

When $\omega_0=3$ and $\Delta=0.1$, by using (\ref{tau}) and (\ref{tau2}) we  draw the Hopf bifurcation values  by  thick black curves   shown in Figure \ref{Tkkkk}(a). The red dashed  curve stands for $-\Delta-2\tau\Delta^2+\frac \tau4k^2=0$,   above which the Hopf bifurcation is subcritical and below which the bifurcation is supercritical. This is  a theoretical proof of FIG.4 in \cite{ar9}.

Four Bautin bifurcation points and two double Hopf points are marked by B$_1-$B$_4$ and HH$_1-$HH$_2$. The blue, thin curves stands for the saddle-node bifurcation of periodic solutions originating from Bautin points by using DDE-Biftool  \cite{ddebif1,ddebif2}.

In fact,  at B$_1$, some calculations yield $k =0.9146,\tau = 0.5288$,  Re$l_1(0,0)=0$ and Re$l_2(0,0)=  -0.1619<0$. At
  B$_2$, we have $k = 0.5819,\tau = 1.5469$,  Re$l_1(0,0)=0$ and Re$l_2(0,0)=-0.0679<0$. Thus from Corollary \ref{corol} and Figure \ref{bautinfig} we know there exists a region near each Bautin points (see Figure\ref{Tkkkk}(b,d)), where a stable equilibrium, a stable periodic orbit and an unstable periodic orbit coexist (Regions VIII or X).

\subsection{Simulations near the double Hopf bifurcation}

Fixing $\omega_0=3,\Delta=0.1$, at the double Hopf point HH$_1$, we have $k = 3.1286$, $\tau = 1.0472$, and two imaginary roots of the characteristic equation (\ref{CE}) are  $\rmi\alpha=-\rmi1.4389$, $\rmi\beta=-\rmi4.5611$, which means the nonresonant condition (\ref{nonre}) is fulfilled.
By using the normal form method given in Section 4, we have $c_1=0.2289 \epsilon -0.6042 \delta$, $c_2=1.9153 \delta +0.2289 \epsilon$,
$d_0=1$, $b_0=c_0=-1.0905$. This corresponds to the  case IVb given in Section 7.5 of \cite{homs}, and there exists two Neimark-Sacker bifurcation curves of periodic solutions (torus bifurcation) originating from HH$_1$.  For sufficiently small $\epsilon$ and $\delta$, the two Neimark-Sacker bifurcation curves can be calculated locally, by $c_2=c_0c_1$ and $c_2=c_1/b_0$, which are $\epsilon=-2.6253 \delta$ and $\epsilon=-3.1018 \delta$.

Using the DDE-Biftool, we draw the global Neimark-Sacker bifurcation curves   in Figure \ref{Tkkkk}(b) in red dash-dotted curves, above which system (\ref{model1}) exhibits unstable oscillations on 2-torus. In Figure \ref{Tkkkk}(c), we give a complete local (near the origin) bifurcation sets in polar coordinates $r_1$ and $r_2$ of Eq.(\ref{NF23polarreducehh}) \cite{homs}. Clearly, in regions III, IV and V, there exist solutions oscillating around both $r_1$-axis and $r_2$-axis.

So far, we have studied the delay coupled Kuramoto oscillators on the OA manifold from the point of view of Bautin and double Hopf bifurcation analysis. Through some elaborative bifurcation analysis, we give a complete bifurcation set in the plane of delay and coupling strength.
To sum up, we can give a result in Table \ref{table1}, where the stability of  $r=0$, number of stable periodic orbits, unstable periodic orbits  or 2-torus of Eq.(\ref{model1}) are shown in the regions I-XI on $\tau-k$ plane. With respect to the Kuramoto model, we recall that $r=0$ represents the incoherence, the periodic solutions stands for the coherent states. Thus we obtain clear scenarios about the synchronous transition   of delayed Kuramoto model (\ref{model}).

\begin{table}[htbp]
\caption{\label{table1} Stability of $r=0$ (s for stable, u for unstable), number of stable periodic orbits, unstable periodic orbits and 2-torus of Eq.(\ref{model1})   for parameters in I-XI.}\begin{center}
\begin{tabular}{@{}lccccccccccc}
\hline
~ & I & II & III & IV & V & VI & VII & VIII & IX & X & XI \\
\hline
 $r=0$  & s & u & u & u & u & u & u & s & s & s & u \\
   stable PO & 2 & 2 & 2 & 2 & 2 & 2 & 1 & 1 & 0 &  1 & 1  \\
   unstable PO & 2 & 1 & 1 & 0 & 1 & 1 & 0 & 1 & 0 & 1 & 0 \\
   2-torus & 0 & 0 & 1 & 1 & 1 & 0 & 0 & 0 & 0 & 0 & 0 \\
\hline
\end{tabular}\end{center}
\end{table}

\subsection{Simulations of the Kuramoto model}

    Now, some illustrations by integrating the Kuramoto model (\ref{model}) are given as  comparisons with the reduced model (\ref{model1}).

\begin{figure}[htbp]
  \centering \includegraphics[width=0.48\textwidth]{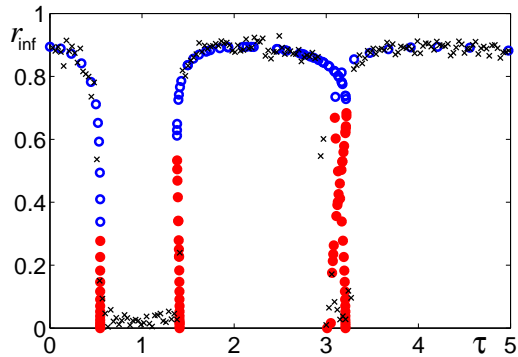}
    \caption{ (color online) $\omega_0=3$, $\Delta=0.1$, and $k=1$. $r_{\textrm{inf}}$ stands for $|r(\infty)|$. The bifurcation diagram of Eq.(\ref{model1})  is shown. Stable periodic orbits are labeled by blue circles, and unstable periodic orbits, red dots.   The order parameters by integrating the Kuramoto model (\ref{model}) with $N=300$ are marked by black crosses. }\label{diagram}
\end{figure}

    Choosing $\omega_0=3$, $\Delta=0.1$, $k=1$, and letting $\tau$ vary, we have a bifurcation diagram shown in Figure \ref{diagram}. By using DDE-Biftool and computing the numbers of Floquet exponents with positive real part, the order parameters $r_{\mathrm{inf}}=|r(\infty)|$ and stability of periodic solutions of (\ref{model1}) (i.e., stability of the coherent states of (\ref{model})) are shown. We find the simulation results by integrating (\ref{model}) coincide with theoretical results of Eq.(\ref{model1}) on the OA manifold very well, and there exist synchrony windows  with the increasing of $\tau$:
    comparing this with Figure \ref{Tkkkk}(b), one can see along these points the transition ``stable periodic solution'' $\rightarrow{\textrm{saddle~node~bifurcation}}\rightarrow$ ``unstable periodic solution'' $\rightarrow{\textrm{subcritical~Hopf~bifurcation}}\rightarrow$ ``stable $r=0$'' $\cdots$.

\begin{figure}[htbp]
  \centering a)\includegraphics[width=0.47\textwidth]{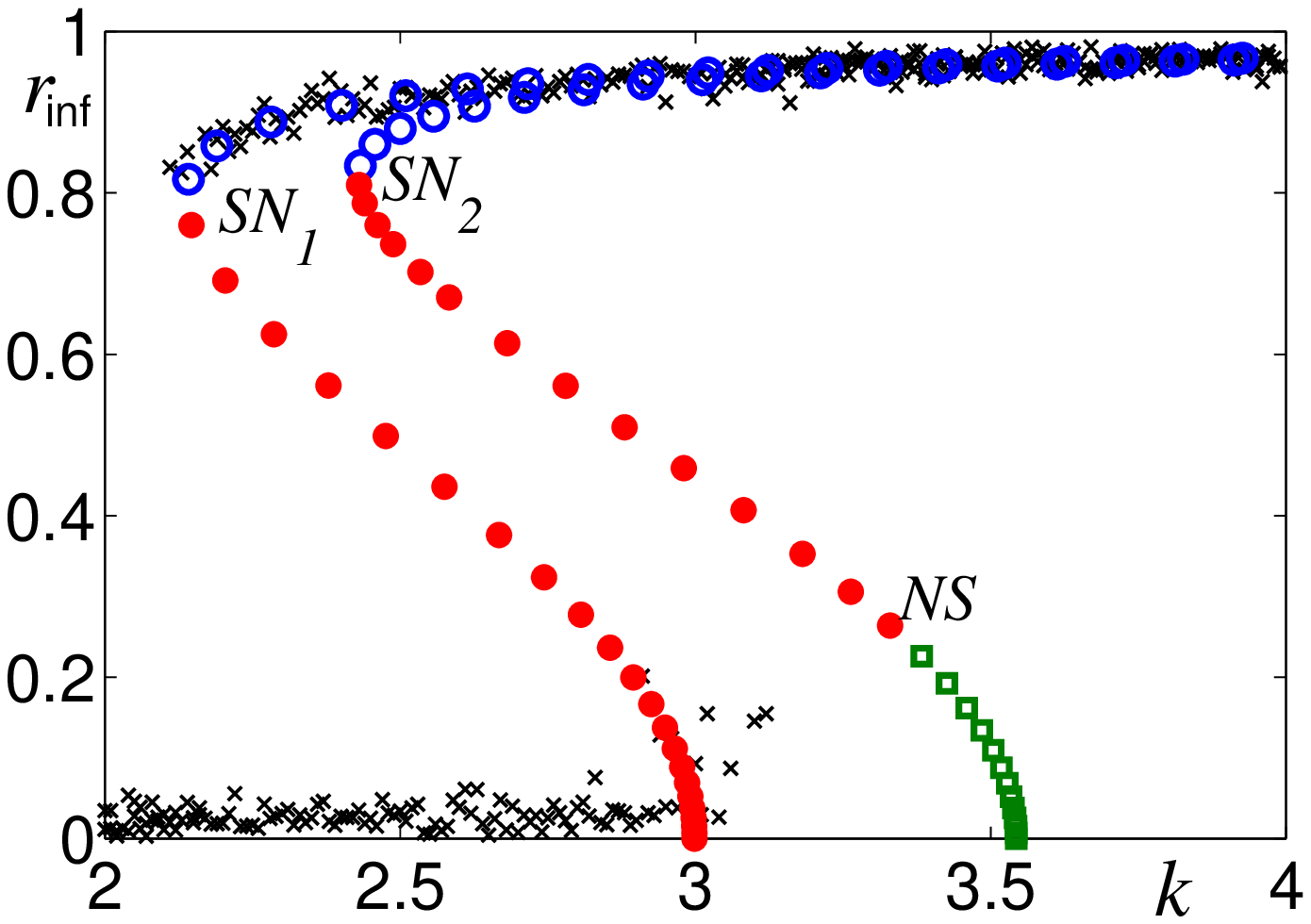} b)\includegraphics[width=0.47\textwidth]{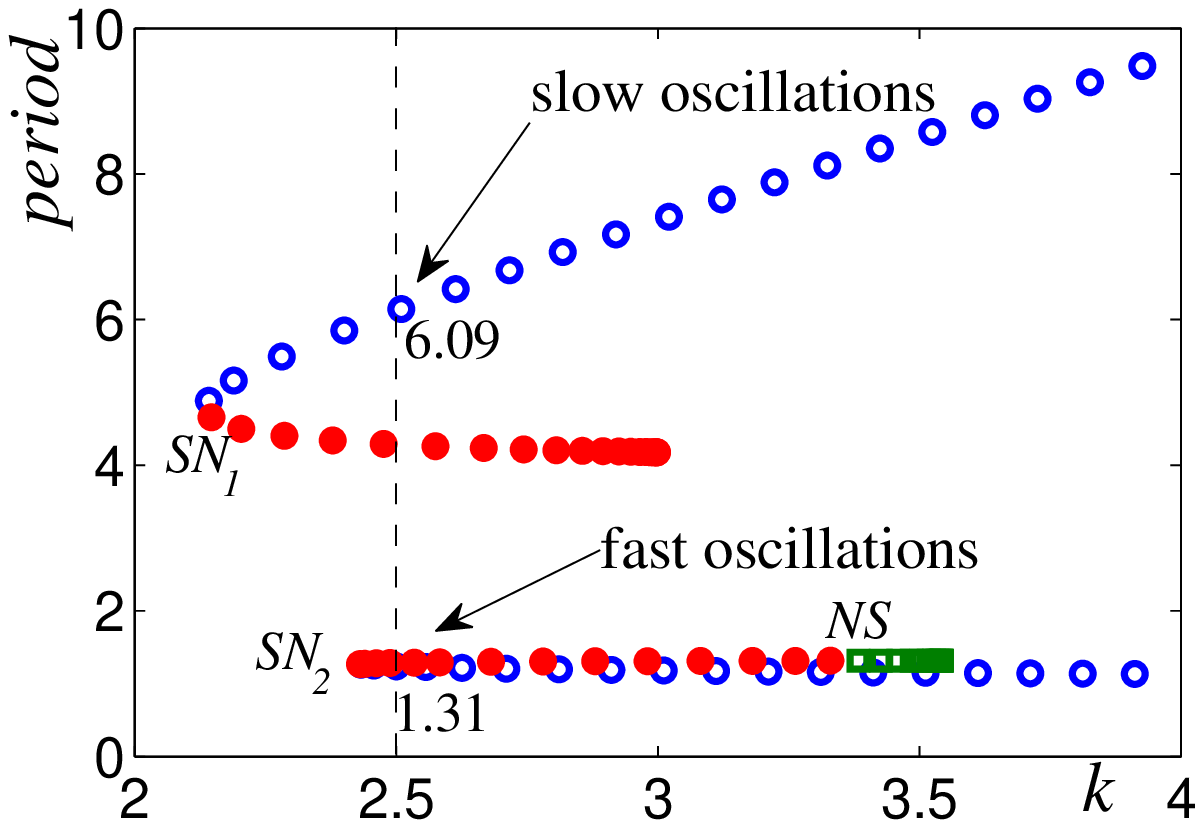}
  \caption{ (color online) $\omega_0=3$, $\Delta=0.1$, and $\tau=1$.  The bifurcation diagram of Eq.(\ref{model1}) and the period of the bifurcating solutions are  shown in a) and b), respectively. $NS$ stands for  Neimark-Sacker  bifurcation. $SN_1$ and $SN_2$ are saddle-node bifurcation of periodic orbits. Stable periodic orbits are labeled by blue circles, unstable (1  Floquet exponent with positive real part) periodic orbits, red dots, and unstable (3  Floquet exponents with positive real part) periodic orbits green squares.   The order parameters by integrating the Kuramoto model (\ref{model}) with $N=300$ are marked by black crosses. }\label{diagram2}
\end{figure}

\begin{figure}[htbp]
  \centering\includegraphics[width=0.45\textwidth]{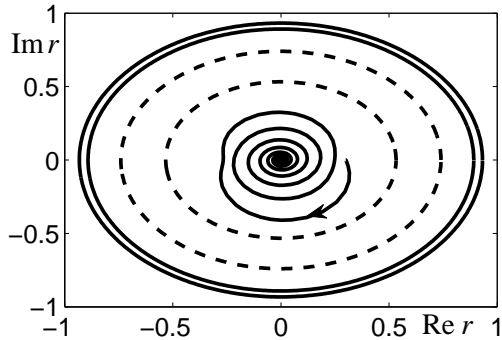}
  \caption{When $k=2.5$ and $\tau=1$, two stable periodic solutions, two unstable periodic solutions (dashed) and  a stable equilibrium  coexist in Eq.(\ref{model1}).}\label{coex}
\end{figure}

\begin{figure}[htbp]
  \centering a)\includegraphics[width=0.28\textwidth]{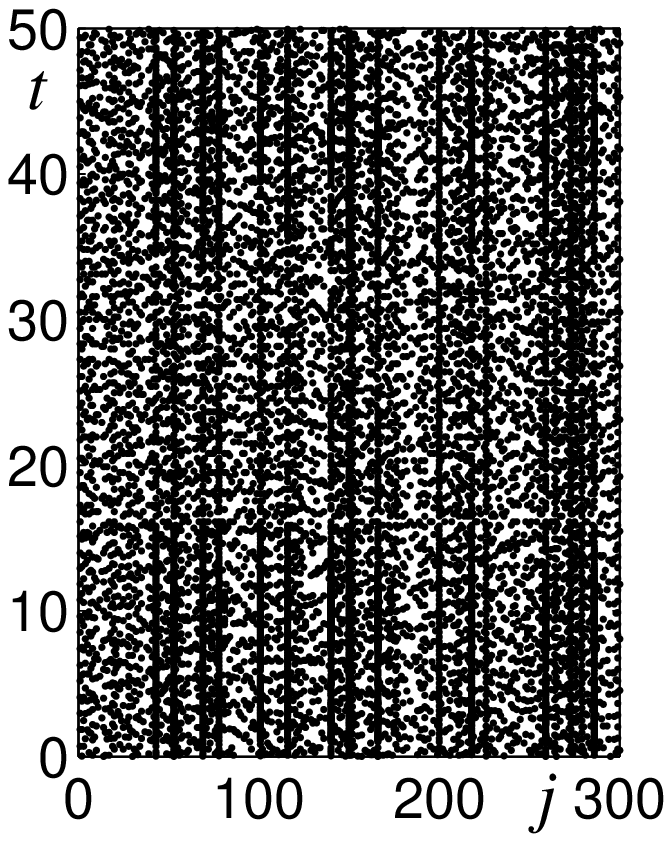}
  b)\includegraphics[width=0.28\textwidth]{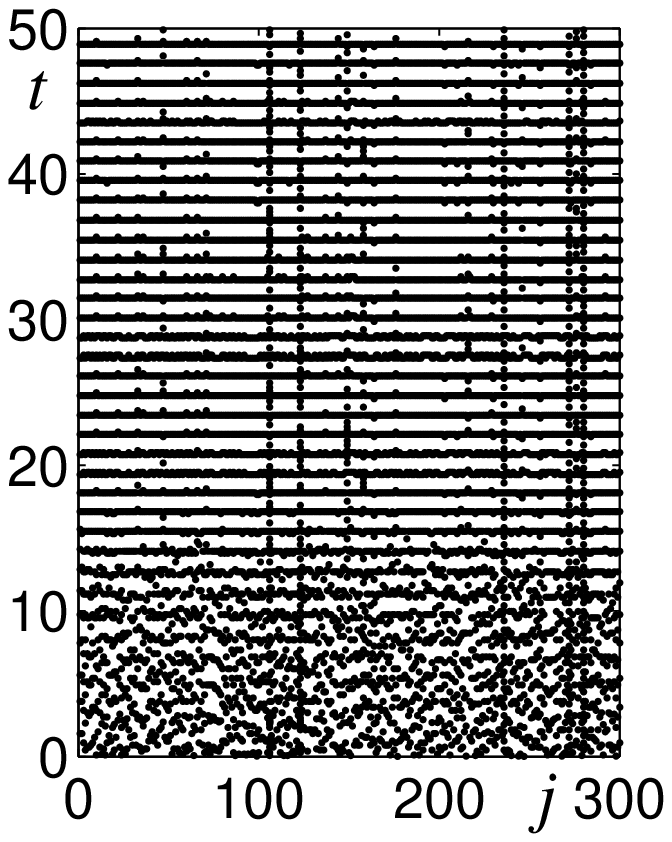}
  c)\includegraphics[width=0.28\textwidth]{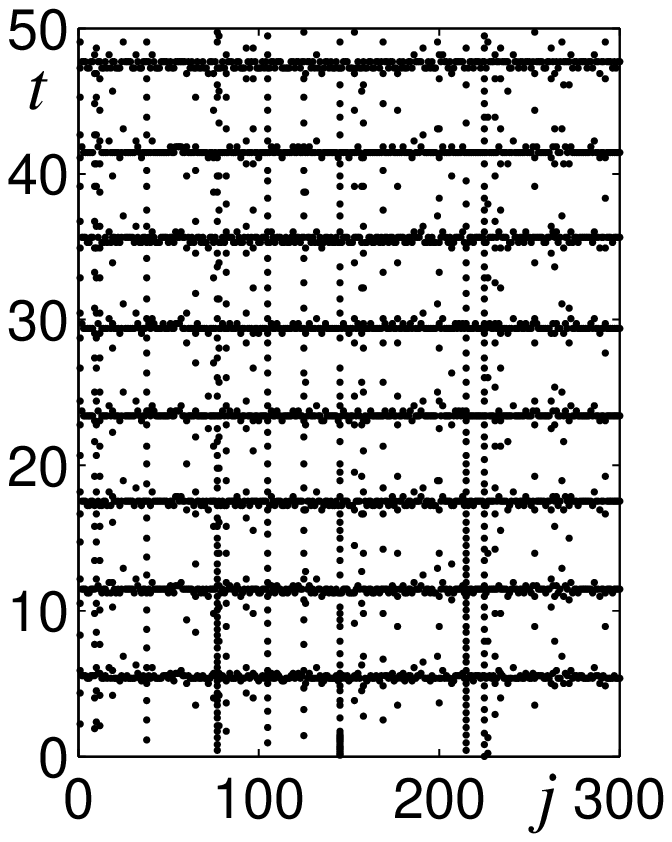}
  \caption{When $k=2.5$ and  $\tau=1$, system (\ref{model}) may exhibit incoherence (with initial values uniformly  distributed on $[0,2\pi]$), coherence of fast oscillation (with initial values uniformly  distributed on $[0,\pi]$), or coherence of slow oscillation (with initial values uniformly  distributed on $[0,0.2\pi]$). In the figures, dots stand for the  points with $\theta_j(t)=2\pi$.}\label{kura}
\end{figure}

     If we fix $\tau=1$ and let $k$ vary, we obtain simulation results in Figure  \ref{diagram2}(a). By using DDE-Biftool and computing the numbers of Floquet exponents with positive real part, the order parameters, together with two saddle-node points and a Neimark-Sacker bifurcation point, are shown. As $k$ increases, the two branches of bifurcating solutions have nearly the same order parameters, but they can be distinguished by periods.  In Figure  \ref{diagram2}(b) we calculate the periods of the two branches of bifurcating periodic solutions, including the fast and slow oscillations.

     In Figure  \ref{diagram2}(a), one can also find from Figure \ref{Tkkkk}(b) that the simulation results coincide with the theoretical results. Between the saddle-node point and the Hopf bifurcation point, stable incoherent and coherent states coexist, i.e., the hysteresis loop. Particularly, we find an interesting phenomenon near $k=2.5$ (in region I), that is  after two times of saddle-node bifurcations,  two stable coherent states, two unstable coherent states and a stable incoherent state coexist, i.e., two hysteresis loops intersect. These coexisting solutions of  (\ref{model1}) are shown by DDE-Biftool in Figure \ref{coex}. Recall that (\ref{model1}) is an infinite-dimensional functional differential equation \cite{Halefde}, thus the unstable periodic
orbits do not separate the two stable periodic orbits.

Fixing $\tau=1$ and $k=2.5$, we perform three simulations about Kuramoto model (\ref{model}) as shown in Figure \ref{kura}, where we find different initial values lead the Kuramoto model to incoherence and coherence respectively.
The two kinds of stable periodic oscillations, shown in Figure \ref{diagram2} and Figure \ref{coex} with periods 1.31 and 6.09,  are simulated in Figure \ref{kura} (b-c).

\subsection{Delay-coupled Hindmarsh-Rose neurons}
Consider the following coupled Hindmarsh-Rose system \cite{HRref,HRref2}
\begin{equation}\label{HR}
\begin{array}{l}
\dot x_j=y_j-x_j^3+3x_j^2-z_j+3\\
\dot y_j=1-5x_j^2-y_j~~~~~~~~~~~~~~~~~~~~~~~~~~~~j=1,2,\ldots,N\\
\dot z_j=0.006[4(x_j+I_j)-z_j]+\frac k N\sum\limits_{i=1}^{N}[z_i(t-\tau)-z_j]
\end{array}
\end{equation}
where the random values $I_j$ are Gaussian distributed with mean 1.56 and variance 0.5. Thus this is a near-identical, delay-coupled system. 

\begin{figure}[htbp]
  \centering a)\includegraphics[width=0.28\textwidth]{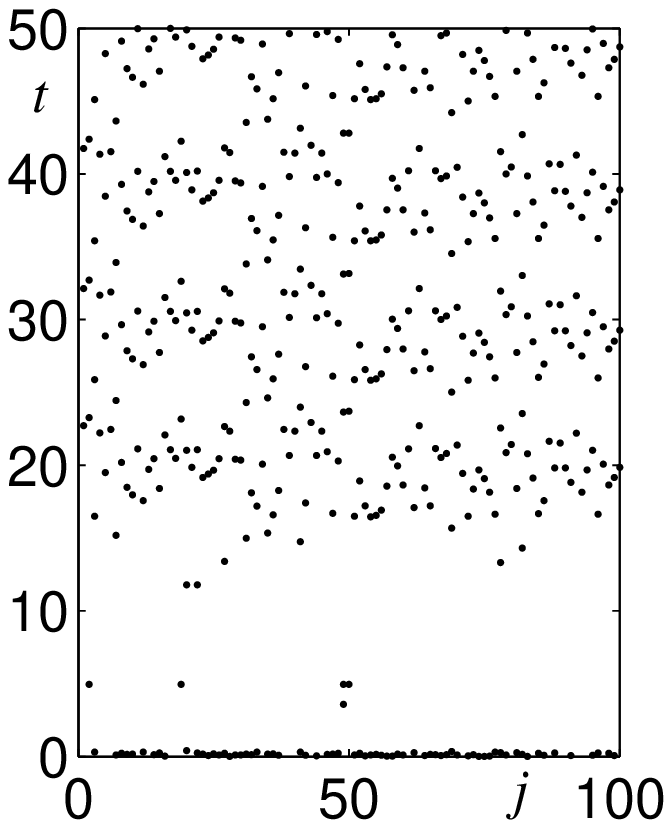}
  b)\includegraphics[width=0.28\textwidth]{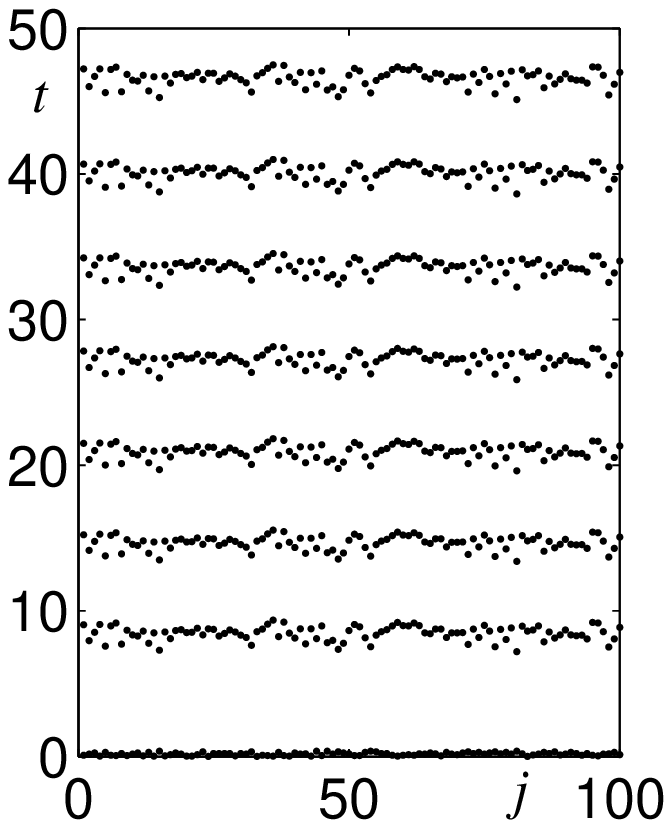}
  c)\includegraphics[width=0.28\textwidth]{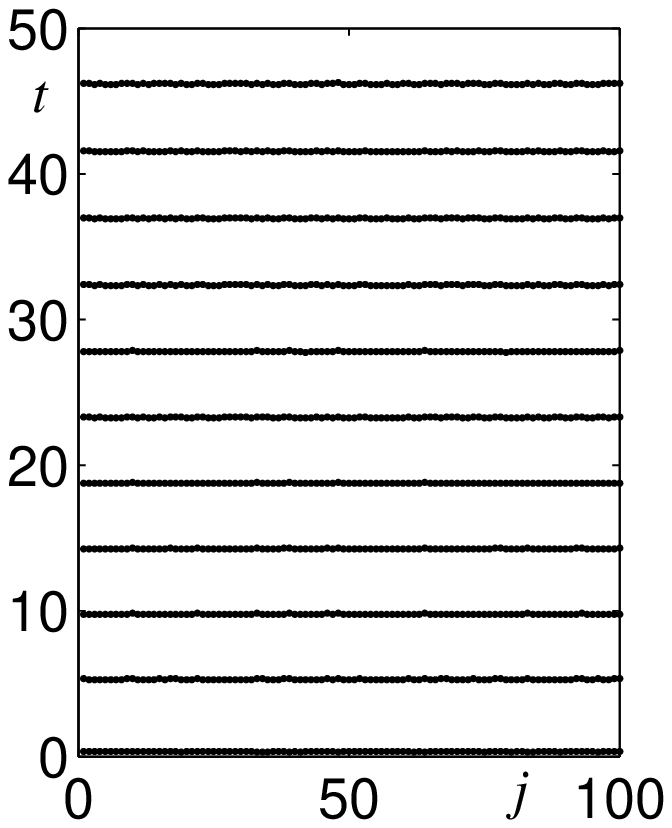}
  \caption{$k=2$,  $\tau=4$ and $t\in[0,50]$. a) When the initial values are uniformly distributed in $[1,4.5]$,  system (\ref{HR}) is incoherence. b) When the initial values are uniformly distributed in $[1,2.5]$, system (\ref{HR}) exhibits  synchronized states with period near 7. c) When the initial values are uniformly distributed in $[1,1.2]$, system (\ref{HR}) exhibits  synchronized states with period near 4.8.  In the figures, dots stand for local maxima of $z_j(t)$.}\label{HRfig}
\end{figure}

Choosing $N=100$, $k=2$ and $\tau=4$, we perform three groups of simulations as shown in Figure \ref{HRfig}.  In the figures we find different initial values  lead the system to the incoherent state or one of at least two coherent states with different periods as well.

\section{Conclusion and Discussion}
In this paper, we study the delay coupled Kuramoto oscillators on the OA manifold from the point of view of Bautin and double Hopf bifurcation analysis. Complete bifurcation sets are given in Figure \ref{Tkkkk} and the existence of coherent and incoherent states are listed in Table \ref{table1}. Through the bifurcation analysis, we find, even in a model as simple as (\ref{model}), the dynamical behavior is complicated: both theoretical investigations and numerical simulations indicate that  Bautin bifurcation and double Hopf bifurcation are very common and must appear in this model, which bring hysteresis loop, multistability and oscillations on torus, respectively. We have theoretically proved that system (\ref{model1}) must undergo  Bautin bifurcations thus hysteresis loop always exists in the delay Kuramoto model (\ref{model}). Particularly, two hysteresis loops may intersect in certain region, which yields that four coherent states (two of which are stable) and a stable incoherence coexist in  Kuramoto model. The coexistence of stable coherent and incoherent states  are all simulated.

Mathematically, an interesting relation is revealed from Corollary \ref{corol}: for small $\tau$, Bautin bifurcation value is $k\propto\tau^{-1/2}$ for  $\Delta\ll1$ and $k\propto\Delta$ for   $\Delta\gg1$. In case of Kuramoto model with near identical oscillators, the former applies. Hence we know that the Bautin bifurcation points moves downwards in the $\tau-k$ plane as the increasing of $\tau$.

 In the current model, we have theoretically proved that, at all Bautin bifurcation points, one always have Re$l_2(0,0)<0$, thus the dynamical behavior near these points are clear as shown in Figure \ref{bautinfig}. However, near the double Hopf points, the situation is not completely clear, because there are twelve kinds of unfoldings near a double Hopf point \cite{homs} and we cannot theoretically determine the sign of $b_0$, $c_0$ $d_0$ and $d_0-b_0c_0$ in general. The numerical example in Section 5 is a special and ``simple'' case of double Hopf bifurcation, but we still find two coexisting synchronized  states in the model. Sometimes, stable 2-torus or 3-torus may appear  after several times of  Neimark-Sacker bifurcations of periodic solutions (e.g., case VIa \cite{homs}), which may  correspond to rotating waves in the Kuramoto model. Even though such situations are complicated, the normal form (\ref{NF23polarreducehh}) can still provide clear bifurcation sets near the double Hopf points.

For practical usage, we simulate a system of delay-coupled Hindmarsh-Rose neurons. We find the results in this paper agree well with a system of delay-coupled Hindmarsh-Rose neurons. In the simulation, we also find two  coexisting stable coherent states and one stable incoherent state. It is worth mentioning that there are recently many results   about the periodic oscillation or phase synchronization arisen from  mathematical biology   \cite{newnewnew,nnn1}.  Linking the results in the current paper and these biological models will be an interesting work and is left as a future study.

\section*{Acknowledgments}
The author deeply appreciates the time and effort that the editor and referees spend on reviewing the manuscript.  
This research is supported by  National Natural Science Foundation of China (11301117 and 11371112), by Heilongjiang Provincial Natural Science Foundation  (QC2014C003) and by the Scientific Research Foundation of Harbin Institute of Technology at Weihai HIT(WH) 201421 and HIT.NSRIF.2016079.

\end{document}